# Numerical Solution of two dimensional coupled viscous Burgers' Equation using the Modified Cubic B-Spline Differential Quadrature Method

*H. S. Shukla*[1], *Mohammad Tamsir*[1*], *Vineet K. Srivastava*[2], *Jai Kumar*[3]

[1]Department of Mathematics & Statistics, DDU Gorakhpur University, Gorakhpur-273009, India
[2]ISRO Telemetry, Tracking and Command Network (ISTRAC), Bangalore-560058, India
[3]ISRO Satellite Center (ISAC), Bangalore-560017, India

**ABSTRACT**

In this paper, a numerical solution of the two dimensional nonlinear coupled viscous Burgers' equation is discussed with the appropriate initial and boundary conditions using the modified cubic B-spline differential quadrature method (MCB-DQM). In this method, the weighting coefficients are computed using the modified cubic B-spline as a basis function in the differential quadrature method. Thus, the coupled Burgers' equations are reduced into a system of ordinary differential equations (ODEs). An optimal five stage and fourth-order strong stability preserving Runge–Kutta (SSP-RK54) scheme is applied to solve the resulting system of ODEs. The accuracy of the scheme is illustrated via two numerical examples. Computed results are compared with the exact solutions and other results available in the literature. Numerical results show that the MCB-DQM is efficient and reliable scheme for solving the two dimensional coupled Burgers' equation.

**Keywords:** Burgers' equation; Reyonlds number; Modified cubic B-spline function; MCB-DQM; SSP-RK54.

## 1. Introduction

Consider the two dimensional nonlinear unsteady coupled viscous Burgers' equations:

$$\frac{\partial u}{\partial t}+u\frac{\partial u}{\partial x}+v\frac{\partial u}{\partial y}=\frac{1}{\mathrm{Re}}\left(\frac{\partial^2 u}{\partial x^2}+\frac{\partial^2 u}{\partial y^2}\right), \tag{1.1}$$

$$\frac{\partial v}{\partial t}+u\frac{\partial v}{\partial x}+v\frac{\partial v}{\partial y}=\frac{1}{\mathrm{Re}}\left(\frac{\partial^2 v}{\partial x^2}+\frac{\partial^2 v}{\partial y^2}\right), \tag{1.2}$$

with the initial conditions:

$$\left.\begin{array}{l} u(x,y,0)=\psi_1(x,y);\ (x,y)\in\Omega \\ v(x,y,0)=\psi_2(x,y);\ (x,y)\in\Omega \end{array}\right\},\tag{1.3}$$

and Dirichlet boundary conditions:

$$\left.\begin{array}{l} u(x,y,t)=\xi(x,y,t);\ (x,y)\in\partial\Omega \\ v(x,y,t)=\zeta(x,y,t);\ (x,y)\in\partial\Omega \end{array}\right\},t>0\tag{1.4}$$

where $\Omega=\{(x,y):a\le x\le b,\ \ c\le x\le d\}$ is the computational domain and $\partial\Omega$ is its boundary, $u(x,y,t)$ and $v(x,y,t)$ are the velocity components to be determined, $\psi_1,\ \psi_2, \xi$ and $\zeta$ are known functions, $\frac{\partial u}{\partial t}$ is unsteady term, $u\frac{\partial u}{\partial x}$ is the nonlinear convection term, $\frac{1}{\mathrm{Re}}\left(\frac{\partial^2 u}{\partial x^2}+\frac{\partial^2 u}{\partial y^2}\right)$ is the diffusion term, and $\mathrm{Re}$ is the Reynolds number.

Coupled viscous Burgers' equation is a more appropriate form of the Navier-Stokes equations having the exact solutions. It has the same convective and diffusion form as the incompressible Navier-Stokes equations, and is a simple model for understanding of various physical flows and problems, such as hydrodynamic turbulence, shock wave theory, wave processes in thermo-elastic medium, vorticity transport, dispersion in porous media[1-3]. Numerical solution of Burgers' equation is a natural and first step towards developing methods for the computation of complex flows. Thus, it has become customary to test new approaches in computational fluid dynamics by implementing novel and new approaches to the Burgers' equation yielding in various finite-differences, finite volume, finite-element and boundary element methods etc.

Analytic solution of two dimensional coupled Burgers' equations was first given by Fletcher[4] using the Hopf-Cole transformation. The numerical solution of coupled Burgers' equations are numerically solved by many researchers[5-15]. In recent years, various researchers[16-22] proposed variant of differential quadrature method for the numerical solution of the one and two dimensional linear/nonlinear differential equations. Korkmaz & Dag[23, 24] proposed cubic B-spline and sinc differential quadrature methods. Arora & Singh[25] proposed the modified cubic B-spline differential quadrature method (MCB-DQM) and applied on one dimensional Burgers' equation to checked its efficiency and accuracy. They found that MCB-DQM is very powerful and efficient scheme as compared to other existing numerical methods. Recently, an extention of the MCB-DQM is proposed by Jiwari & Yuan[26] to show the computational modeling of two-dimensional reaction–diffusion Brusselator model with appropriate initial and Neumann boundary conditions. Researchers[27-29] developed an optimal strong stability preserving (SSP) high order time discretization schemes. Strong stability properties of SSP methods is preserved in any norm, semi norm or convex functional of the spatial discretization coupled with the first order Euler time stepping. A description of the optimal explicit and implicit SSP Runge-Kutta and multistep methods is also discussed by the authors.

In this paper, we study the numerical simulation of the two-dimensional unsteady nonlinear coupled viscous Burgers' equations for different Reynolds number. The efficacy and accuracy of the method is confirmed by taking two test problem with suitable initial and boundary conditions. This study shows that the MCB-DQM results are acceptable and in good agreement with the exact solutions and earlier results available in the literature. Rest of the article is prepared as: In Section 2, the modified cubic B-spline differential quadrature method is introduced. In Section 3, the implementation procedure for the problem (1.1) – (1.2) together with the initial conditions (1.3) and boundary conditions (1.4) is illustrated; In Section 4, two test problems are given to establish the applicability and accuracy of the method, while the Section 5 concludes our study.

## 2. Modified cubic B-spline differential quadrature method

In 1972, Bellman et al.[16] introduced differential quadrature method (DQM). This method approximates the spatial derivatives of a function using the weighted sum of the functional values at the certain discrete points. In DQM, the weighting coefficients are determined using several kinds of test functions such as spline function[23], sinc function[24], Lagrange interpolation polynomials, Legendre polynomials[17-22] etc. This section revisits the MCB-DQM[25-26] in order to complete our problem in two dimension. It is assumed that the $M$ and $N$ grid points: $a = x_1 < x_2, ..... < x_M = b$ and $c = y_1 < y_2, ..... < y_N = d$ are uniformly distributed with the spatial step size $\Delta x = x_{i+1} - x_i$ and $\Delta y = y_{j+1} - y_j$ in $x$ and $y$ directions, respectively.

The first and second order spatial partial derivatives of $u(x, y, t)$ with respect to $x$ (keeping $y_j$ as fixed) and with respect to $y$ (keeping $x_i$ as fixed), approximated at $x_i$ and $y_j$, respectively, are defined as:

$$\frac{\partial u(x_i, y_j, t)}{\partial x} = \sum_{k=1}^{N} w_{ik}^{(1)} u(x_k, y_j, t), \qquad i = 1, 2, ..., M \tag{2.1}$$

$$\frac{\partial^2 u\left(x_i, y_j, t\right)}{\partial x^2} = \sum_{k=1}^{N} w_{ik}^{(2)} u\left(x_k, y_j, t\right), \qquad i = 1, 2, ..., M \tag{2.2}$$

$$\frac{\partial u\left(x_i, y_j, t\right)}{\partial y} = \sum_{k=1}^{M} \overline{w}_{jk}^{(1)} u\left(x_i, y_k, t\right), \qquad j = 1, 2, ..., N \tag{2.3}$$

$$\frac{\partial^2 u\left(x_i, y_j, t\right)}{\partial y^2} = \sum_{k=1}^{M} \overline{w}_{jk}^{(2)} u\left(x_i, y_k, t\right), \qquad j = 1, 2, ..., N. \tag{2.4}$$

In the same way, the first and second order spatial partial derivatives of $v(x, y, t)$ with respect to $x$ and with respect to $y$ are approximated as:

$$\frac{\partial v\left(x_i, y_j, t\right)}{\partial x} = \sum_{k=1}^{N} w_{ik}^{(1)} v\left(x_k, y_j, t\right), \qquad i = 1, 2, ..., M \tag{2.5}$$

$$\frac{\partial^2 v\left(x_i, y_j, t\right)}{\partial x^2} = \sum_{k=1}^{N} w_{ik}^{(2)} v\left(x_k, y_j, t\right), \qquad i = 1, 2, ..., M \tag{2.6}$$

$$\frac{\partial v\left(x_i, y_j, t\right)}{\partial y} = \sum_{k=1}^{M} \overline{w}_{jk}^{(1)} v\left(x_i, y_k, t\right), \qquad j = 1, 2, ..., N \tag{2.7}$$

$$\frac{\partial^2 v\left(x_i, y_j, t\right)}{\partial y^2} = \sum_{k=1}^{M} \overline{w}_{jk}^{(2)} v\left(x_i, y_k, t\right), \qquad j = 1, 2, ..., N, \tag{2.8}$$

where $w_{ij}^{(r)}$ and $\overline{w}_{ij}^{(r)}$, $r = 1, 2$ are the weighting coefficients of the *r*th-order spatial partial derivatives with respect to $x$ and $y$.

The cubic B-spline basis functions[22] at the knots are defined as:

$$\varphi_m\left(x\right) = \frac{1}{h^3}\begin{cases} \left(x - x_{m-2}\right)^3, & x \in \left(x_{m-2}, x_{m-1}\right) \\ \left(x - x_{m-2}\right)^3 - 4\left(x - x_{m-1}\right)^3, & x \in \left(x_{m-1}, x_m\right) \\ \left(x_{m+2} - x\right)^3, & x \in \left(x_{m+1}, x_{m+2}\right) \\ 0, & \text{otherwise,} \end{cases} \tag{2.9}$$

where $\{\varphi_0, \varphi_1, ..., \varphi_N, \varphi_{N+1}\}$ is chosen in such a way that it forms a basis over the domain $\Omega = \{(x, y):\ a \le x \le b;\ c \le y \le d\}$. The values of cubic B-splines and its derivatives at the nodal points are depicted in **Table 1**.

**Table 1:** Coefficients of the cubic B-spline $\varphi_m$ and its derivatives at the node $x_m$.

| | $x_{m-2}$ | $x_{m-1}$ | $x_m$ | $x_{m+1}$ | $x_{m+2}$ |
|---|---|---|---|---|---|
| $\varphi_m(x)$ | 0 | 1 | 4 | 1 | 0 |
| $\varphi_m'(x)$ | 0 | $3/h$ | 0 | $-3/h$ | 0 |
| $\varphi_m''(x)$ | 0 | $6/h^2$ | $-12/h^2$ | $6/h^2$ | 0 |

Now, to get a diagonally dominant system of the linear equations, the cubic B-spline basis functions are modified as[25]:

$$\left.\begin{aligned} \phi_1(x) &= \varphi_1(x) + 2\varphi_0(x) \\ \phi_2(x) &= \varphi_2(x) - \varphi_0(x) \\ \phi_m(x) &= \varphi_m(x) \; for \; \mathrm{m} = 3,\dots,N-2 \\ \phi_{N-1}(x) &= \varphi_{N-1}(x) - \varphi_{N+1}(x) \\ \phi_N(x) &= \varphi_N(x) + 2\varphi_{N+1}(x) \end{aligned}\right\}, \tag{2.10}$$

where $\{\phi_1, \phi_2, \dots, \phi_N\}$ forms a basis over the domain $\Omega = \{(x, y):\; a \le x \le b;\; c \le y \le d\}$.

In Eq. (2.1), substituting the values of $\phi_m(x), \quad m = 1, 2, \dots, N$, we get a system of linear equations:

$$\phi_m'(x_i) = \sum_{k=1}^{N} w_{ik}^{(1)} \phi_m(x_k), \; \mathrm{i} = 1, 2, \dots, M. \tag{2.11}$$

With the help of Eq. (2.9), (2.10) and **Table 1**, Eq. (2.11) reduces into a tridiagonal system of equations:

$$A\vec{w}^{(1)}[i] = \vec{R}[i], \; for \; i = 1, 2, \dots, M, \tag{2.12}$$

where $\vec{w}^{(1)}[i] = \left[w_{i1}^{(1)}, w_{i2}^{(1)}, \dots, w_{iN}^{(1)}\right]^T$ is the weighting coefficient vector corresponding to $x_i$, $\vec{R}[i] = \left[\phi_{1,i}', \phi_{2,i}', \dots \dots, \phi_{N-1,i}', \phi_{N,i}'\right]^T$, and the coefficient matrix $A$ is given by:

$$A = \begin{bmatrix} \phi_{1,1} & \phi_{1,2} & & & & & \\ \phi_{2,1} & \phi_{2,2} & \phi_{2,3} & & & & \\ & \phi_{3,2} & \phi_{3,3} & \phi_{3,4} & & & \\ & & \ddots & \ddots & \ddots & & \\ & & & \phi_{N-2,N-3} & \phi_{N-2,N-2} & \phi_{N-2,N-1} & \\ & & & & \phi_{N-1,N-2} & \phi_{N-1,N-1} & \phi_{N-1,N} \\ & & & & & \phi_{N,N-1} & \phi_{N,N} \end{bmatrix}$$

Here, we point out that the coefficient matrix $A$ is invertible. The tridiagonal system of linear equations (2.12) is solved for each $i$ using Thomas algorithm, which gives the weighting coefficients $w_{ik}^{(1)}$ of the first order partial derivative. In a similar way, the weighting coefficients $w_{ij}^{(2)}, 1 \le i, j \le N$ are determined. Weighting coefficients $w_{ij}^{(2)}, 1 \le i, j \le N,$ can be computed as[22]:

$$\begin{cases} w_{ij}^{(r)} = r\left( w_{ij}^{(1)} w_{ii}^{(r-1)} - \dfrac{w_{ij}^{(r-1)}}{x_i - x_j} \right), \textit{ for } i \ne j \text{ and } i = 1,2,3,...,N; \quad r = 2,3,...,N-1 \\ w_{ii}^{(r)} = -\sum\limits_{j=1, j\ne i}^{N} w_{ij}^{(r)}, \text{ for } i = j, \end{cases} \tag{2.13}$$

where $\bar{w}_{ij}^{(r-1)}$ and $\bar{w}_{ij}^{(r)}$ are the weighting coefficients of the $(r-1)^{\text{th}}$ and $r^{\text{th}}$ order partial derivatives with respect to *x*.

Similarly, the weighting coefficients $\bar{w}_{jk}^{(1)}$ of the first order partial derivatives with respect to $y$ using the modified cubic B-Spline functions in Eq. (11) is obtained. Weighting coefficients $\bar{w}_{ij}^{(2)}, 1 \le i, j \le N$ for the second derivatives can be computed from the formule:

$$\begin{cases} \bar{w}_{ij}^{(r)} = r\left( \bar{w}_{ij}^{(1)} \bar{w}_{ii}^{(r-1)} - \dfrac{\bar{w}_{ij}^{(r-1)}}{x_i - x_j} \right), \textit{ for } i \ne j \text{ and } i = 1,2,3,...,N; \quad r = 2,3,...,N-1 \\ \bar{w}_{ii}^{(r)} = -\sum\limits_{j=1, j\ne i}^{N} \bar{w}_{ij}^{(r)}, \text{ for } i = j, \end{cases} \tag{2.14}$$

where $\bar{w}_{ij}^{(r-1)}$ and $\bar{w}_{ij}^{(r)}$ are the weighting coefficients of the $(r-1)^{\text{th}}$ and $r^{\text{th}}$ order partial derivatives with respect to *y*.

## 3. MCB-DQM for two-dimensional coupled Burgers' equation

On substituting the approximate values of the spatial derivatives computed by the MCB-DQM, Eq. (1.1) can be written as:

$$\frac{\partial u(x_i, y_j, t)}{\partial t} = -u(x_i, y_j)\sum_{k=1}^{M} w_{ik}^{(1)} u(x_k, y_j) - v(x_i, y_j)\sum_{k=1}^{N} \bar{w}_{jk}^{(1)} u(x_i, y_k)$$
$$+\frac{1}{\text{Re}}\left[\sum_{k=1}^{M} w_{ik}^{(2)} u(x_k, y_j) + \sum_{k=1}^{N} \bar{w}_{jk}^{(2)} u(x_i, y_k)\right], (x_i, y_j) \in R,\ t > 0,\ i = 1,2,...,M, \quad j = 1,2,...,N. \tag{3.1}$$

Similarly, Eq. (2) can be written as:

$$\frac{\partial u\left(x_i,y_j,t\right)}{\partial t}=-u\left(x_i,y_j\right)\sum_{k=1}^{M}w_{ik}^{(1)}v\left(x_k,y_j\right)-v\left(x_i,y_j\right)\sum_{k=1}^{N}\overline{w}_{jk}^{(1)}v\left(x_i,y_k\right)$$
$$+\frac{1}{\mathrm{Re}}\left[\sum_{k=1}^{M}w_{ik}^{(2)}v\left(x_k,y_j\right)+\sum_{k=1}^{N}\overline{w}_{jk}^{(2)}v\left(x_i,y_k\right)\right],\left(x_i,y_j\right)\in R,\ t>0,\ i=1,2,...,M,\quad j=1,2,...,N. \tag{3.2}$$

Eq. (3.1) and Eq. (3.2) reduce into a system of ordinary differential equations:

$$\frac{du\left(x_i,y_j,t\right)}{dt}=F_1\left(u\left(x_i,y_j,t\right)\right),\ \ i=1,2,...,M \text{ and } j=1,2,...,N. \tag{3.3}$$

$$\frac{dv\left(x_i,y_j,t\right)}{dt}=F_2\left(u\left(x_i,y_j,t\right)\right),\ \ i=1,2,...,M \text{ and } j=1,2,...,N. \tag{3.4}$$

Eqs. (3.3) and (3.4) together with the initial conditions (1.3) and Dirichlet boundary conditions (1.4) are solved by SSP-RK54 scheme.

## 4. Results and discussion

Here we consider two test problems of two dimensional coupled Burgers equation as given in the introduction part to provide the MCB-DQM numerical solutions. The accuracy and consistency of the scheme is measured in terms of error norms $L_2$ and $\mathrm{L}_\infty$, defined as:

$$\left.\begin{aligned}&L_2:=\| u_{exact}-u_{computed}\|_2=\sqrt{\sum_{i=0}^{n}\sum_{j=0}^{n}| u_{i,j}^{\ exact}-u_{i,j}^{\ computed}|^2}\\&\mathrm{L}_\infty:=\| u_{exact}-u_{computed}\|_\infty=\max_{i,j}| u_{i,j}^{\ exact}-u_{i,j}^{\ computed}|\end{aligned}\right\}, \tag{4.1}$$

where $u_{exact}$ and $u_{computed}$ represent the exact and computed solutions, respectively.

Numerical solutions to Eqs. (1) and (2) will be tested for the following two test problems.

**4.1. Problem:** The analytical solutions of Eqs. (1.1) and (1.2) can be generated as[4]:

$$\left.\begin{aligned}&u\left(x,y,t\right)=\frac{3}{4}-\frac{1}{4\left[1+\exp\left(\left(-4x+4y-t\right)\mathrm{Re}/32\right)\right]}\\&v\left(x,y,t\right)=\frac{3}{4}+\frac{1}{4\left[1+\exp\left(\left(-4x+4y-t\right)\mathrm{Re}/32\right)\right]}\end{aligned}\right\}\ \left(\mathrm{x,y}\right)\in\Omega, \tag{4.2}$$

The square domain $\Omega = \{(x,y): 0 \le x \le 1, 0 \le y \le 1\}$ is considered as the computational domain, and the initial and boundary conditions for $u(x,y,t)$ and $v(x,y,t)$ are taken from the analytical solutions Eq. (4.2).

**4.2. Problem:** In this problem, we take the computational domain

$$\Omega = \{(x,y): 0 \le x \le 0.5, 0 \le y \le 0.5\}$$

with the initial conditions:

$$\left.\begin{aligned} u(x,y,0) &= \sin(\pi x) + \cos(\pi y) \\ v(x,y,0) &= x + y \end{aligned}\right\} \quad (\mathrm{x,y}) \in \Omega, \tag{4.3}$$

and Dirichet boundary conditions:

$$\left.\begin{aligned} u(0,y,t) &= \cos(\pi y), \\ u(0.5,y,t) &= 1 + \cos(\pi y), \\ v(0,y,t) &= y, \\ v(0.5,y,t) &= 0.5 + y, \end{aligned}\right\} \; 0 \le y \le 0.5,\; t \ge 0, \tag{4.4}$$

$$\left.\begin{aligned} u(x,0,t) &= 1 + \sin(\pi x), \\ u(x,0.5,t) &= \sin(\pi x), \\ v(x,0,t) &= x, \\ v(x,0.5,t) &= x + 0.5, \end{aligned}\right\} \; 0 \le x \le 0.5,\; t \ge 0, \tag{4.5}$$

For the test problem 4.1, we have taken a grid size $20 \times 20$ with time step $\Delta t = 0.0001$ and $\mathrm{Re} = 100$. Computed and exact values of $u$ and $v$ are shown in Tables 2 and 3 along with the results given by Srivastava et al.[14-15] and Bahadir[9] at some typical grid point. The tabulated results show that the proposed scheme produces better result than Bahadir[9]. Tables 4 and 5 show the errors $L_2$ and $\mathrm{L}_\infty$, and also the rate of convergence of $u$ and $v$ components, respectively, at $\mathrm{Re} = 100$, $t = 1.0$ for $\Delta t = 0.0001$. From Tables 4 and 5, it can be observed that the MCB-DQM performs better than Srivastava et al.[15] and gives more than quadratic rate of convergence. The MCB-DQM computed solutions of $u$ and $v$ for $\mathrm{Re} = 100$ at $t = 0.5$ are depicted in Fig. 1 while Fig. 2 shows analytical solutions of $u$ and $v$, respectively.

For the test problem 4.2, numerical computations are carried out with the parameters: $\mathrm{Re} = 50$, 100, and $20 \times 20$ grid, $t = 0.625$ for $\Delta t = 0.0001$ in order to compare the computed results with those given by Jain & Holla[5], Bahadir[9], and Srivastava et al.[14-15]. Table 6 shows the comparisons of numerical results obtained using the MCB-DQM scheme at $t = 0.625$, with the methods of Jain & Holla[5], Bahadir[9], and Srivastava et al.[14-15]. From Table 6, it can be noticed that the computed MCB-DQM results are in good agreement with Jain & Holla[5], Bahadir[9], and Srivastava et al.[14-15]. Fig. 3 depicts the MCB-DQM computed $u$ and $v$ solutions corresponding to $\mathrm{Re}$ =50, 100 and 500 at $t = 0.625$.

**Table 2:** Comparison of MCB-DQM and exact solutions of $u$ for $\mathrm{Re} = 100$, $20 \times 20$ grid and $\Delta t = 0.0001$.

| Grid $(x, y)$ | $t = 0.5$ | | | | | $t = 2.0$ | | | | |
|---|---|---|---|---|---|---|---|---|---|---|
| | MCB-DQM | Exact | I-LFDM[14] | Expo-FDM[15] | Bahadir[9] | MCB-DQM | Exact | I-LFDM[14] | Expo-FDM[15] | Bahadir[9] |
| (0.1,0.1) | 0.54412 | 0.54332 | 0.54300 | 0.54300 | 0.54235 | 0.50050 | 0.50048 | 0.50047 | 0.50047 | 0.49983 |
| (0.5,0.1) | 0.50037 | 0.50035 | 0.50034 | 0.50034 | 0.49964 | 0.50000 | 0.50000 | 0.50000 | 0.50000 | 0.49930 |
| (0.9,0.1) | 0.50000 | 0.50000 | 0.50000 | 0.50000 | 0.49931 | 0.50000 | 0.50000 | 0.50000 | 0.50000 | 0.49930 |
| (0.3, 0.3) | 0.54388 | 0.54332 | 0.54269 | 0.54270 | 0.54207 | 0.50050 | 0.50048 | 0.50044 | 0.50044 | 0.49977 |
| (0.7, 0.3) | 0.50037 | 0.50035 | 0.50032 | 0.50032 | 0.49961 | 0.50000 | 0.50000 | 0.50000 | 0.50000 | 0.49930 |
| (0.1, 0.5) | 0.74196 | 0.74221 | 0.74215 | 0.74215 | 0.74130 | 0.55632 | 0.55568 | 0.55515 | 0.55516 | 0.55461 |
| (0.5, 0.5) | 0.54347 | 0.54332 | 0.54251 | 0.54252 | 0.54222 | 0.50050 | 0.50048 | 0.50041 | 0.50041 | 0.49973 |
| (0.9, 0.5) | 0.50035 | 0.50035 | 0.50030 | 0.50030 | 0.49997 | 0.50001 | 0.50000 | 0.50000 | 0.50000 | 0.49931 |
| (0.3, 0.7) | 0.74211 | 0.74221 | 0.74211 | 0.74212 | 0.74146 | 0.55597 | 0.55568 | 0.55482 | 0.55482 | 0.55429 |
| (0.7, 0.7) | 0.54327 | 0.54332 | 0.54246 | 0.54247 | 0.54243 | 0.50054 | 0.50048 | 0.50038 | 0.50038 | 0.49970 |
| (0.1, 0.9) | 0.74994 | 0.74995 | 0.74994 | 0.74994 | 0.74913 | 0.74406 | 0.74426 | 0.74420 | 0.74420 | 0.74340 |
| (0.5, 0.9) | 0.74219 | 0.74221 | 0.74210 | 0.74210 | 0.74201 | 0.55575 | 0.55568 | 0.55450 | 0.55451 | 0.55413 |
| (0.9, 0.9) | 0.54333 | 0.54332 | 0.54228 | 0.54229 | 0.54232 | 0.50052 | 0.50048 | 0.50053 | 0.50053 | 0.50001 |

**Table 3:** Comparison of MCB-DQM and exact solutions of $v$ for $\mathrm{Re} = 100$, $20 \times 20$ grid and $\Delta t = 0.0001$.

| Grid $(x, y)$ | $t = 0.5$ | | | | | $t = 2.0$ | | | | |
|---|---|---|---|---|---|---|---|---|---|---|
| | MCB-DQM | Exact | I-LFDM[14] | Expo-FDM[15] | Bahadir[9] | MCB-DQM | Exact | I-LFDM[14] | Expo-FDM[15] | Bahadir[9] |
| (0.1,0.1) | 0.95589 | 0.95668 | 0.95700 | 0.95700 | 0.95577 | 0.99950 | 0.99952 | 0.99953 | 0.99953 | 0.99826 |
| (0.5,0.1) | 0. 99963 | 0.99965 | 0.99966 | 0.99966 | 0.99827 | 1.00000 | 1.00000 | 1.00000 | 1.00000 | 0.99860 |
| (0.9,0.1) | 1.00000 | 1.00000 | 1.00000 | 1.00000 | 0.99861 | 1.00000 | 1.00000 | 1.00000 | 1.00000 | 0.99861 |
| (0.3, 0.3) | 0.95612 | 0.95668 | 0.95731 | 0.95731 | 0.95596 | 0.99950 | 0.99952 | 0.99956 | 0.99956 | 0.99820 |
| (0.7, 0.3) | 0.99964 | 0.99965 | 0.99968 | 0.99968 | 0.99827 | 1.00000 | 1.00000 | 1.00000 | 1.00000 | 0.99860 |
| (0.1, 0.5) | 0.75804 | 0.75779 | 0.75785 | 0.75785 | 0.75699 | 0.94368 | 0.94432 | 0.94485 | 0.94485 | 0.94393 |
| (0.5, 0.5) | 0.95654 | 0.95668 | 0.95749 | 0.95749 | 0.95685 | 0.99950 | 0.99952 | 0.99959 | 0.99959 | 0.99821 |
| (0.9, 0.5) | 0.99965 | 0.99965 | 0.99970 | 0.99970 | 0.99903 | 0.99999 | 1.00000 | 1.00000 | 1.00000 | 0.99862 |
| (0.3, 0.7) | 0.75789 | 0.75779 | 0.75789 | 0.75789 | 0.75723 | 0.94403 | 0.94432 | 0.94518 | 0.94518 | 0.94409 |
| (0.7, 0.7) | 0.95673 | 0.95668 | 0.95754 | 0.95754 | 0.95746 | 0.99946 | 0.99952 | 0.99962 | 0.99962 | 0.99823 |
| (0.1, 0.9) | 0.75006 | 0.75005 | 0.75006 | 0.75006 | 0.74924 | 0.75595 | 0.75574 | 0.75580 | 0.75580 | 0.75500 |
| (0.5, 0.9) | 0.75781 | 0.75779 | 0.75790 | 0.75790 | 0.75781 | 0.94425 | 0.94432 | 0.94550 | 0.94550 | 0.94441 |
| (0.9, 0.9) | 0.95667 | 0.95668 | 0.95772 | 0.95772 | 0.95777 | 0.99948 | 0.99952 | 0.99948 | 0.99948 | 0.99846 |

**Table 4.** $L_2$, $L_\infty$ errors and rate of convergence for the $u$-component for $\mathrm{Re}=100$, $\Delta t=$ 0.0001 at $t=1.0$.

| Grid | $L_2$ | | | $L_\infty$ | | |
|---|---|---|---|---|---|---|
| | Expo-FDM[15] | MCB-DQM | | Expo-FDM[15] | MCB-DQM | |
| | | | Rate | | | Rate |
| $4\times4$ | 8.5708e-002 | 1.6388e-002 | - | 9.7046e-002 | 2.8788e-003 | |
| $8\times8$ | 4.9429e-002 | 1.9286e-003 | 3.0875 | 4.6886e-002 | 1.9572e-004 | 3.8786 |
| $16\times16$ | 1.9192e-002 | 3.9474e-004 | 2.2881 | 2.0467e-002 | 2.0486e-005 | 3.2561 |
| $32\times32$ | 8.6812e-003 | 8.1181e-005 | 2.2817 | 9.0744e-003 | 2.2202e-006 | 3.2059 |
| $64\times64$ | - | 1.5322e-005 | 2.4055 | - | 2.1838e-007 | 3.3458 |

**Table 5**. $L_2$, $L_\infty$ errors and rate of convergence for the $v$-component for $\mathrm{Re}=100$, $\Delta t=$ 0.0001 at $t=1.0$.

| Grid | $L_2$ | | | $L_\infty$ | | |
|---|---|---|---|---|---|---|
| | Expo-FDM[15] | MCB-DQM | | Expo-FDM[15] | MCB-DQM | |
| | | | Rate | | | Rate |
| $4\times4$ | 8.5708e-002 | 1.6388e-002 | - | 9.7046e-002 | 2.8788e-003 | - |
| $8\times8$ | 4.9431e-002 | 1.9286e-003 | 3.0875 | 4.6887e-002 | 1.9573e-004 | 3.8786 |
| $16\times16$ | 1.9196e-002 | 3.9474e-004 | 2.2881 | 2.0471e-002 | 2.0486e-005 | 3.2561 |
| $32\times32$ | 8.6878e-003 | 8.1181e-005 | 2.2817 | 9.0813e-003 | 2.2202e-006 | 3.2059 |
| $64\times64$ | - | 1.5322e-005 | 2.4055 | - | 2.1838e-007 | 3.3458 |

**Table 6:** Comparison of computed results of $u$ and $v$ for $\mathrm{Re}=50$, grid size $20\times20$ and $\Delta t=$ 0.0001 at $t=0.625$.

| Grid $(x, y)$ | Computed values of $u$ | | | | | Computed values of $v$ | | | | |
|---|---|---|---|---|---|---|---|---|---|---|
| | MCB-DQM | I-LFDM[14] | Expo-FDM[15] | Bahadir[9] | Jain and Holla[5] | MCB-DQM | I-LFDM[14] | Expo-FDM[15] | Bahadir[9] | Jain and Holla[5] |
| (0.1, 0.1) | 0.97056 | 0.97146 | 0.97146 | 0.96688 | 0.97258 | 0.09842 | 0.09869 | 0.09869 | 0.09824 | 0.09773 |
| (0.3, 0.1) | 1.15152 | 1.15280 | 1.15280 | 1.14827 | 1.16214 | 0.14107 | 0.14158 | 0.14158 | 0.14112 | 0.14039 |
| (0.2, 0.2) | 0.86244 | 0.86308 | 0.86308 | 0.85911 | 0.86281 | 0.16732 | 0.16754 | 0.16754 | 0.16681 | 0.16660 |
| (0.4, 0.2) | 0.98078 | 0.97985 | 0.97985 | 0.97637 | 0.96483 | 0.17223 | 0.17111 | 0.17111 | 0.17065 | 0.17397 |
| (0.1, 0.3) | 0.66336 | 0.66316 | 0.66316 | 0.66019 | 0.66318 | 0.26380 | 0.26378 | 0.26378 | 0.26261 | 0.26294 |
| (0.3, 0.3) | 0.77226 | 0.77233 | 0.77233 | 0.76932 | 0.77030 | 0.22653 | 0.22655 | 0.22655 | 0.22576 | 0.22463 |
| (0.2, 0.4) | 0.58273 | 0.58181 | 0.58181 | 0.57966 | 0.58070 | 0.32935 | 0.32851 | 0.32851 | 0.32745 | 0.32402 |
| (0.4, 0.4) | 0.76179 | 0.75862 | 0.75862 | 0.75678 | 0.74435 | 0.32884 | 0.32502 | 0.32502 | 0.32441 | 0.31822 |

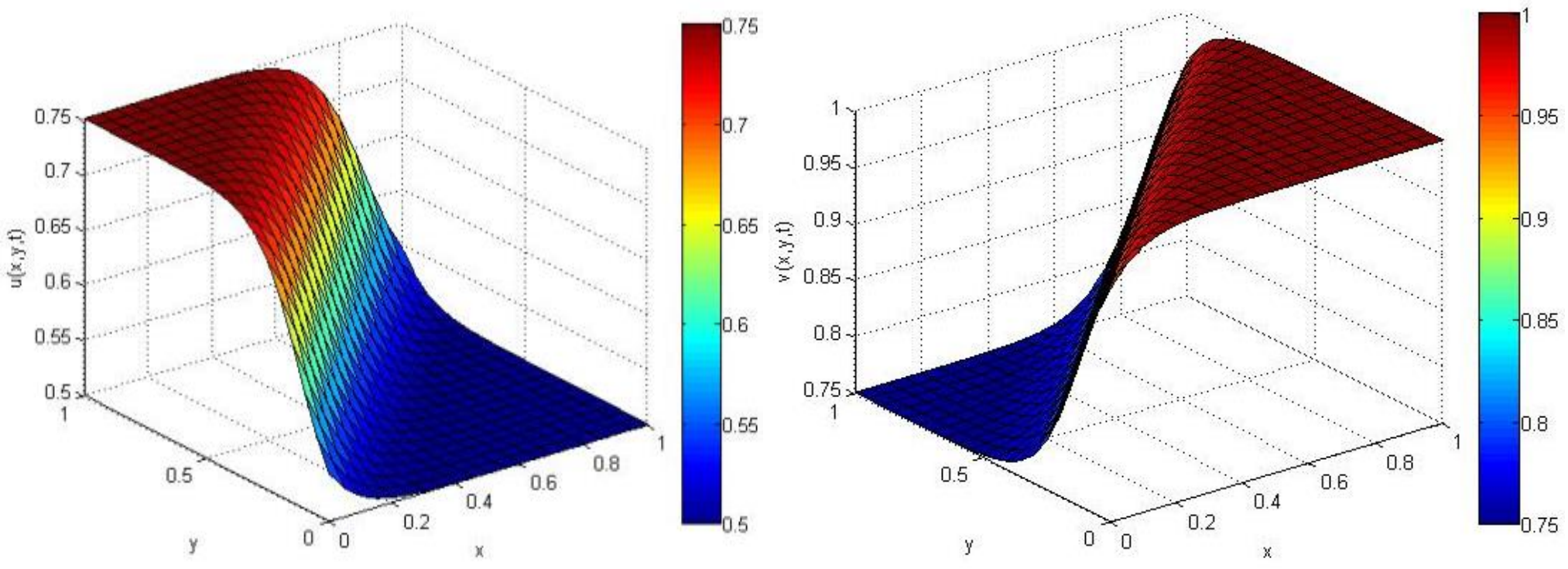

**Fig. 1:** Numerical solution at $t = 0.5$ with $\Delta t = 0.0001$, $\mathrm{Re} = 100$ and grid size $20 \times 20$ for the test problem 4.1.

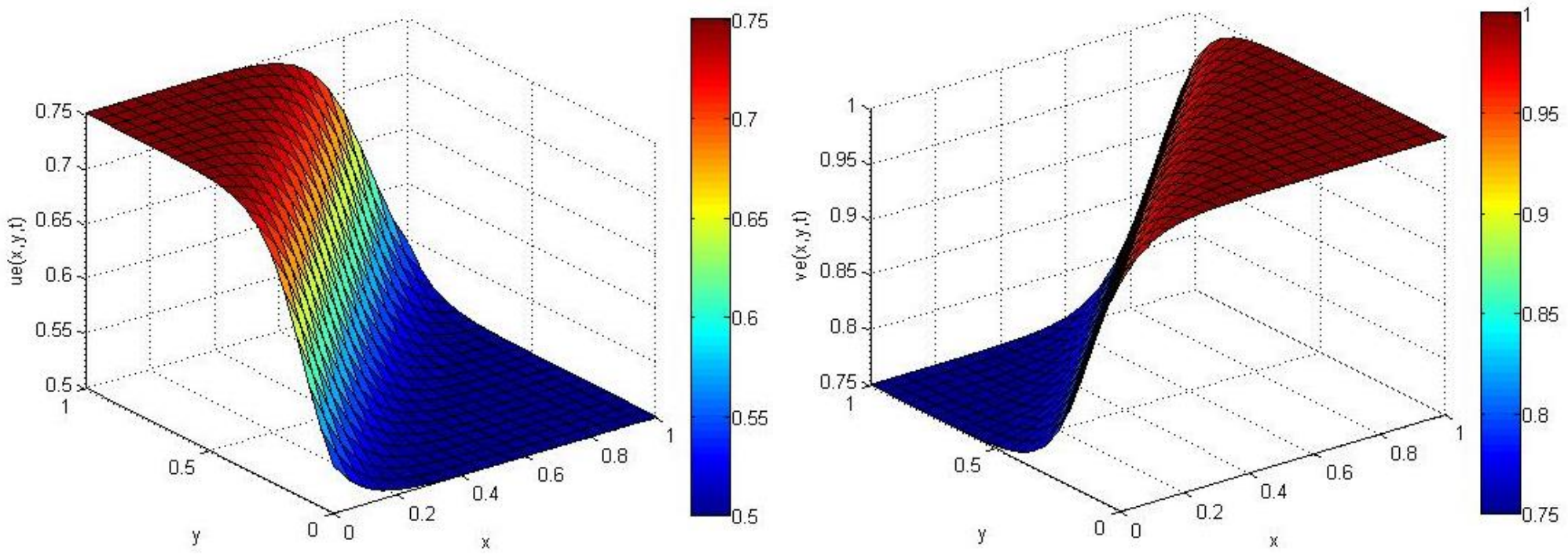

**Fig. 2.** Exact solution at $t = 0.5$ with $\Delta t = 0.0001$, $\mathrm{Re} = 100$ and grid size $20 \times 20$ for the test problem 4.1.

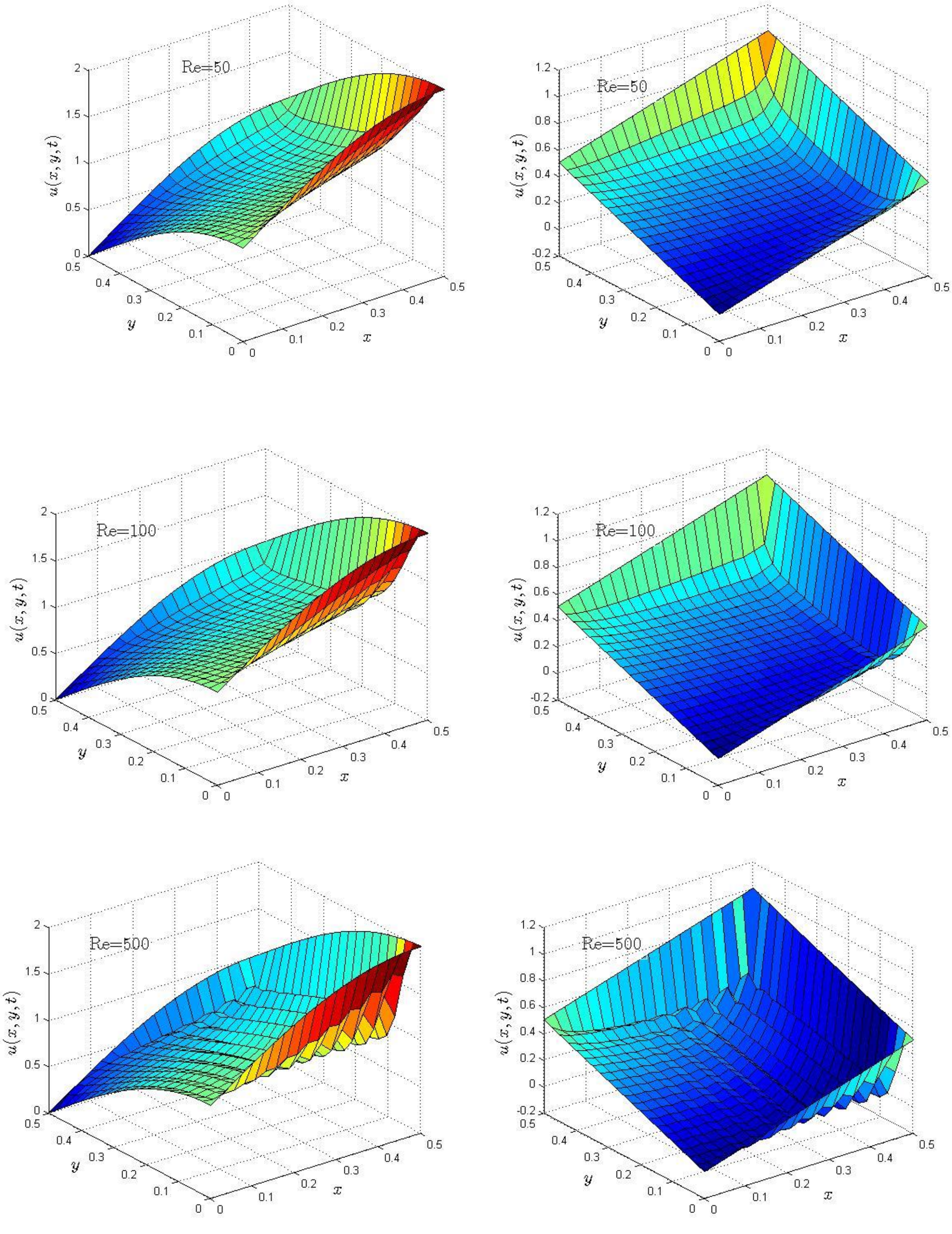

**Fig. 3.** Numerical solutions at $t = 0.625$ with $\Delta t = 0.0001$ and grid size $20 \times 20$ for the problem 4.2.

## 5. Conclusions

A modified cubic B-spline differential quadrature method is presented for the numerical solutions of two dimensional nonlinear coupled viscous Burgers' equations. The computed results show that the solution obtained by this scheme is highly accurate and very close to the exact solutions. We also notice that the scheme has more than quadratic rate of convergance. The obtained results show that the MCB-DQM is a promising numerical scheme for solving the higher dimensional nonlinear physical problems governed by partial differential equations.

## References

[1]Cole J. D., "On a quasilinear parabolic equations occurring in aerodynamics", Quart Appl Math 9, 225 (1951).

[2]Esipov S. E., "Coupled Burgers' equations: a model of poly-dispersive sedimentation", Phys Rev 52, 3711 (1995).

[3]J. D. Logan, "An introduction to nonlinear partial differential equations", Wily-Interscience, New York (1994).

[4]C. A. J. Fletcher, "Generating exact solutions of the two-dimensional Burgers' equation", Int. J. Numer. Meth. Fluids 3, 213 (1983).

[5]P. C. Jain, D. N. Holla, "Numerical solution of coupled Burgers' equations", Int. J. Numer. Meth.Eng. 12, 213 (1978).

[6]C. A. J. Fletcher, "A comparison of finite element and finite difference of the one- and two-dimensional Burgers' equations", J. Comput. Phys. 51, 159 (1983).

[7]F.W. Wubs, E.D. de Goede, "An explicit–implicit method for a class of time-dependent partial differential equations", Appl. Numer. Math. 9,157 (1992).

[8]O. Goyon, "Multilevel schemes for solving unsteady equations", Int. J. Numer. Meth. Fluids 22, 937 (1996).

[9]Bahadir A. R., "A fully implicit finite-difference scheme for two-dimensional Burgers' equation", Applied Mathematics and Computation 137, 131 (2003).

[10]V. K. Srivastava, M. Tamsir, U. Bhardwaj, Y. Sanyasiraju, "Crank-Nicolson scheme for numerical solutions of two dimensional coupled Burgers' equations", IJSER 2(5), 44 (2011).

[11]M. Tamsir, V. K. Srivastava, "A semi-implicit finite-difference approach for two-dimensional coupled Burgers' equations", IJSER 2(6), 46 (2011).

[12]V. K. Srivastava, M. Tamsir, "Crank-Nicolson semi-implicit approach for numerical solutions of two-dimensional coupled nonlinear Burgers' equations", Int. J. Appl. Mech. Eng. 17 (2), 571 (2012).

[13]V. K. Srivastava, M. K. Awasthi, M. Tamsir,"A fully implicit Finite-difference solution to one-dimensional Coupled Nonlinear Burgers' equations", Int. J. Math. Sci. 7(4), 23 (2013).

[14]V. K. Srivastava, M. K. Awasthi, S. Singh,"An implicit logarithm finite difference technique for two dimensional coupled viscous Burgers' equation", AIP Advances 3, 122105 (2013).

[15]V. K. Srivastava, S. Singh, M. K. Awasthi, "Numerical solutions of coupled Burgers' equations by an implicit finite-difference scheme", AIP Advances 3, 082131 (2013).

[16]R. Bellman, B. G. Kashef, J. Casti, Differential quadrature: a technique for the rapid solution of nonlinear differential equations, J. Comput. Phy. 10, 40-52, (1972).

[17]C. Shu, B. E. Richards, Application of generalized differential quadrature to solve two dimensional incompressible navier-Stokes equations, Int. J. Numer. Meth. Fluids, 15, 791-798, (1992).

[18]J. R. Quan, C.T. Chang, New insights in solving distributed system equations by the quadrature methods-I, Comput. Chem. Eng. 13, 779–788, (1989).

[19]J. R. Quan, C.T. Chang, New insights in solving distributed system equations by the quadrature methods-II, Comput. Chem. Eng. 13, 1017–1024, (1989).

[20]C. Shu, Y.T. Chew, Fourier expansion-based differential quadrature and its application to Helmholtz eigenvalue problems, Commun. Numer. Methods Eng. 13 (8), 643–653, (1997).

[21]C. Shu, H. Xue, Explicit computation of weighting coefficients in the harmonic differential quadrature, J. Sound Vib. 204 (3), 549–555, (1997).

[22]C. Shu, Differential Quadrature and its Application in Engineering, Athenaeum Press Ltd., Great Britain, (2000).

[23]A. Korkmaz, I. Dag, Cubic B-spline differential quadrature methods and stability for Burgers' equation, Eng. Comput. Int. J. Comput. Aided Eng. Software, 30 (3), 320–344, (2013).

[24]A. Korkmaz, I. Dag, Shock wave simulations using sinc differential quadrature method, Eng. Comput. Int. J. Comput. Aided Eng. Software, 28(6), 654–674, (2011).

[25]G. Arora, B. K. Singh, Numerical solution of Burgers' equation with modified cubic B-spline differential quadrature method, Applied Math. Comput., 224 (1), 166-177, (2013).

[26]R. Jiwari, J. Yuan, A computational modeling of two dimensional reaction–diffusion Brusselator system arising in chemical processes, J. Math. Chem., 52 (6), 1535-1551, (2014).

[27]S. Gottlieb, C. W. Shu, E. Tadmor, Strong Stability-Preserving High-Order Time Discretization Methods, SIAM REVIEW, 43(1), 89-112, (2001).

[28]J. R. Spiteri, S. J. Ruuth, A new class of optimal high-order strong stability-preserving time-stepping schemes, SIAM J. Numer. Analysis. 40 (2), 469-491, (2002).
[29]S. Gottlieb, D. I. Ketcheson, C. W. Shu, High Order Strong Stability Preserving Time Discretizations, J. Sci. Comput., 38, 251-289, (2009).